\documentclass{article}
\usepackage[T1]{fontenc}
\usepackage[utf8]{inputenc}
\usepackage{microtype}
\usepackage{amsmath,amssymb,amsfonts,amsbsy,mathtools}
\usepackage{enumitem}
\usepackage[numbers,sort&compress]{natbib}
\usepackage[colorlinks = true,linkcolor = blue,urlcolor  = blue,citecolor = blue,anchorcolor = blue]{hyperref}
\numberwithin{equation}{section}
\newtheorem{theorem}{Theorem}[section]
\newtheorem{proposition}[theorem]{Proposition}

\newtheorem{remark}[theorem]{Remark}

\newcommand{\R}{\mathbb{R}}
\newcommand{\E}{\mathbb{E}}
\newcommand{\Pp}{\mathbb{P}}
\newcommand{\F}{\mathcal{F}}
\newcommand{\dd}{\mathrm{d}}

\begin{document}

\title{Conditional Path Decomposition at the Infimum and Maximum Drawdowns for Spectrally Negative L\'{e}vy Processes}

\author{ C. Vardar-Acar\footnote{Middle East Technical University, Department of Statistics,  \"{U}niversiteler Mah. Dumlupınar Blv. No:1, 06800 \c{C}ankaya Ankara, Turkey, \texttt{cvardar@metu.edu.tr}} \and M. \c{C}a\u{g}lar\footnote{ College of Sciences,  Department of Mathematics,
Ko\c{c} University, Rumeli Feneri Yolu, 34450 Sariyer, Istanbul Turkey,
\texttt{mcaglar@ku.edu.tr }} }
%\author{
%\name{Ceren Vardar-Acar\textsuperscript{a}\thanks{Ceren Vardar-Acar Email:cvardar@metu.edu.tr}, Mine \c{C}a\u{g}lar\textsuperscript{b}\thanks{Mine \c{C}a\u{g}lar Email:mcaglar@ku.edu.tr}}
%\affil{\small Department of Statistics, Middle East Technical University, Ankara, T\"urkiye; \textsuperscript{a} Department of Mathematics, Ko\c{c} University, \.{I}stanbul, T\"urkiye. \textsuperscript{b}}
%}
%\authorone[Department of Statistics, Institute of Applied Mathematics, Middle East Technical University, Ankara, T\"urkiye. ]{Ceren Vardar Acar} 
%\authortwo[Department of Mathematics, Ko\c{c} University, \.{I}stanbul, T\"urkiye. ]{Mine \c{C}a\u{g}lar}

%\addressone{Dumlupınar Bulvarı, No:1, 06800 Çankaya, Ankara, Türkiye} % Your postal address goes here.
%\emailone{cvardar@metu.edu.tr} %Authors email goes here.
%\addresstwo{Rumelifeneri Yolu 34450 Sarıyer, İstanbul, Türkiye} % Your postal address goes here.
%\emailtwo{mcaglar@ku.edu.tr} %Authors email goes here.

\maketitle

\begin{abstract}
We study  maximum-drawdown laws conditioned on extremes for a spectrally negative L\'evy process and observed up to an independent exponential time. The main contribution is a set of scale-function characterizations of the pre-infimum path arising from two decompositions of the process. The first is the decomposition at the infimum into pre-infimum and post-infimum components. The second, under the ordering in which the infimum is attained before the supremum, decomposes the path into pre-infimum, intermediate, and post-supremum components. We also identify the distribution of the supremum for the pre-infimum process in the first decomposition. The resulting conditional laws are expressed as Doob $h$-transforms of killed spectrally negative L\'evy processes and they yield explicit formulas for the maximum drawdown on each  independent path component.  The results confirm the classical decompositions for Brownian motion.
\end{abstract}


\section{Introduction}
% Drawdown processes have a long history in risk management and mathematical finance.  Portfolio-selection problems under drawdown constraints were studied by \cite{GrossmanZhou1993} and later generalized by \cite{CvitanicKaratzas1995}.  Drawdown-based risk measures and their use in portfolio optimization appear in \cite{Chekhlov2005}.  The drawdown process is often defined as a risk measure, and using the drawdown process, performance measures such as the Calmar ratio, Sterling ratio, and Burke ratio have been introduced to compare the performance of stocks \cite{schuhmacher2011sufficient}.

% \cite{Pospisil2009} introduced a new Greek variable dependent on the drawdown variable, which was used to measure the sensitivity of an investment. \cite{Carr2011} developed a European-style digital drawdown insurance contract and proposed its use in futures transactions to reduce the exchange rate risk of insurance claims. There are many more studies on drawdown insurance available in the literature such as (\cite{Palmowski2018}, \cite{Landriault2015}, \cite{Landriault2016}, \cite{Azcue2005}, \cite{Carr2011}). The Russian option, constructed based on the supremum values of its underlying asset, was defined in \cite{shepp1993russian}. For this put option, the optimal stopping time $T^*$ is found as the first time the drawdown exceeds a certain level $k^*$.
Maximum drawdown is considered as a measure of the largest realized loss/risk in the financial instruments.   Its conditional distribution is relevant for risk assessment, timing rules, drawdown-linked contracts, and path-dependent derivatives, (\cite{Pospisil2010}, \cite{Hadjiliadis2006}, \cite{Zhang2010}, \cite{Chekhlov2005}, \cite{cvitanic1995portfolio}, \cite{grossman1993optimal}, \cite{Chekhlov2005}, \cite{Magdon-Ismail2004}, \cite{Leal2005}, \cite{Carr2011}, \cite{Vecer2007}, \cite{Huang2013}, \cite{Rossello2021}, \cite{Sornette2003}, \cite{Vecer2006}, \cite{Avram2004}, \cite{Mijatovic2012}, \cite{Landriault2015} and \cite{Baurdoux2017}). The distributional properties of maximum drawdowns are studied in (\cite{Douady2000}, \cite{Magdon-Ismail2004}, \cite{Salminen2007},  \cite{VardarAcar2017}, \cite{VardarAcar2021}, \cite{SALMINEN20205592}, \cite{salminen2025drawdowns}, \cite{ZHANG2023104669}, \cite{10.1214/ECP.v20-3945}, \cite{risks7040105}). 
 The distributions of maximum drawdown are closely connected to first-passage problems. Scale functions provide tractable identities for such passage problems. In the spectrally negative L\'{e}vy setting,  see  \cite{Avram2020}, \cite{Bertoin1996}, \cite{Kyprianou2014} and \cite{Kuznetsov2013}. For standard Brownian motion, \cite{Douady2000} and \cite{Salminen2007} have obtained laws of maximum drawdown and maximum drawup. Explicitly in Brownian and diffusion models, distributions are found up to a passage time as well as exponential time. Maximum drawdown and maximum drawup have been analyzed  for spectrally negative L\'{e}vy processes until the passage time of a given level in \cite{VardarAcar2017}.  Path decompositions at extrema for spectrally negative L\'{e}vy processes together with applications to maximum drawdown and drawdown duration were developed by \cite{VardarAcar2021}. However, the law of the pre-infimum component was left as an open question in a setting where the process was split first into two at its global infimum, and next into three pieces under the condition that infimum is reached  before the supremum. The present paper completes the missing part by deriving the law of the pre-infimum component  and obtaining the distribution of its maximum-drawdown. 
 
In this paper, we consider a spectrally negative L\'{e}vy process $X$  observed up to an independent exponential time $T$ with rate $\gamma>0$. We assume that $X$  is not the negative of a subordinator and hence it is of unbounded variation. Let
\[
S_t:=\sup_{0\le s\le t}X_s,\qquad I_t:=\inf_{0\le s\le t}X_s
\]
denote its running supremum and infimum processes, respectively, for $0\le t\le T$.  The drawdown and drawup processes are defined as
\[
D_t:=S_t-X_t,\qquad U_t:=X_t-I_t,
\]
and the corresponding maximum drawdown and maximum drawup are
\begin{eqnarray}\label{eq:MDD-MDU-def}
M_t^-&:&=\sup_{0\le u\le v\le t}(X_u-X_v)=\sup_{0\le s\le t}D_s,\nonumber\\
M_t^+&:&=\sup_{0\le u\le v\le t}(X_v-X_u)=\sup_{0\le s\le t}U_s.
\end{eqnarray}
Conditional on the  extrema
\[
I_T=a<0<S_T=b,
\]
the path is decomposed according to infimum, and also with respect to both extremes given that the  infimum is attained before the  supremum.  Explicitly, we focus on the ordering
\[
H_I<H_S,
\]
where $H_I$ and $H_S$ denote   attainment  times of $I_T$ and $S_T$, respectively.  In this case, the path decomposes into three pieces: the pre-infimum, the intermediate path from the infimum to the supremum, and the post-supremum.  

The novelty of the present work is the identification of the conditional law of the pre-infimum path component towards the objective of obtaining the conditional distribution of the maximum drawdown. We characterize this component both in the decomposition at the infimum and also, more delicately, under the  conditioning $H_I < H_S, I_T=a, S_T=b,$
where the infimum is attained before the supremum. These conditional pre-infimum laws are obtained using Doob-$h$ transforms and in terms of well known scale functions $W^{(\gamma)}$, $Z^{(\gamma)}$. The formulas include the derivative of  $W^{(\gamma)}$, which exists everywhere by our assumption that $X$ is of unbounded  variation. Explicitly,  we show in Theorem \ref{Prop_pre_post_inf_dist} that 
the law of the pre-$H_I$ process is an $h$-transform of the law of the spectrally negative L\'{e}vy process killed at  $\tau_a\wedge T $ with 
\[
h(x) = \frac{\gamma}{\Phi(\gamma)} W^{(\gamma)\prime}(x-a) -\gamma W^{(\gamma)}(x-a)
\] 
for $x\ge a$, obtain the conditional distribution of the supremum  as
	\begin{eqnarray*}
		% &&\mathbb{P}_0\left(S_{H_{I}} < b \mid I_{T}=a\right)=\frac{ \frac{\gamma}{\Phi(\gamma)} W^{(\gamma)\prime}(0) }{\frac{\gamma}{\Phi(\gamma)} W^{(\gamma)\prime}(-a) -\gamma W^{(\gamma)}(-a)}\left[Z^{(\gamma)}(-a) - Z^{(\gamma)}(b-a) \frac{W^{(\gamma)}(-a)}{W^{(\gamma)}(b-a)}\right]
        &&\mathbb{P}_0\left\{S_{H_{I}} < b \mid I_{T}=a\right\}=1-\frac{ \frac{\gamma}{\Phi(\gamma)} W^{(\gamma)\prime}(b-a) -\gamma W^{(\gamma)}(b-a) }{\frac{\gamma}{\Phi(\gamma)} W^{(\gamma)\prime}(-a) -\gamma W^{(\gamma)}(-a)}\left[\frac{W^{(\gamma)}(-a)}{W^{(\gamma)}(b-a)}\right]
	\end{eqnarray*}
	for  $b\ge 0$, and derive the conditional distribution of  maximum drawdown $M_{0, H_{I}}^{-}$ of the pre-$H_I$ process given by
\begin{eqnarray*}
\Pp\left\{M_{0, H_{I}}^{-}< d \mid I_{T}=a \right\} =1-\frac{ \frac{\gamma}{\Phi(\gamma)} W^{(\gamma)\prime}(d) -\gamma W^{(\gamma)}(d) }{\frac{\gamma}{\Phi(\gamma)} W^{(\gamma)\prime}(-a) -\gamma W^{(\gamma)}(-a)}\left[\frac{W^{(\gamma)}(-a)}{W^{(\gamma)}(d)}\right]
\end{eqnarray*}
for $d\ge -a$. In Proposition \ref{prop:HI-before-HS}, the process is split at   both infumum, with $I_T=a$, and  supremum, with $S_T=b$. In this case, the pre-$H_I$ distribution is found as 
 the Doob-$h$ transform of the process killed on leaving $(a,b)$ through the lower boundary, with
\begin{equation*}\label{eq:h1-HIHS}
h_1(x)=
\frac{W^{(\gamma)\prime}(x-a)W^{(\gamma)}(b-a)-W^{(\gamma)\prime}(b-a)W^{(\gamma)}(x-a)}
     {\big(W^{(\gamma)}(b-a)\big)^2},
\qquad a<x<b.
\end{equation*}
Under the same decomposition, the distribution of maximum drawdown  is found in Theorem \ref{theorem:mdd-HIHS}  as
\begin{equation*}\hspace{-1cm}\small\label{eq:f1-piecewise}
\Pp_{a,b}^{<}\{M_{0,H_I}^-<d\}=\begin{cases}
0, & d\le -a,\\
\displaystyle
1-\left[\frac{W^{(\gamma)\prime}(d)W^{(\gamma)}(r)-W^{(\gamma)\prime}(r)W^{(\gamma)}(d)}{W^{(\gamma)\prime}(-a)W^{(\gamma)}(r)-W^{(\gamma)\prime}(r)W^{(\gamma)}(-a) }\right]
\left(\frac{W^{(\gamma)}(-a)}{W^{(\gamma)}(d)}\right),
& -a<d<r,\\
1, & d\ge r.
\end{cases}
\end{equation*}
where $r=b-a$.   

The paper is organized as follows. The preliminary definitions and our assumptions are given in Section \ref{sec:preliminaries}. We perform two path decompositions through the extremes, one at the global infimum in Section \ref{sec:infimum-decomposition} and one at both infimum and supremum in Section \ref{sec:ordered-decomposition}. In Theorem  \ref{Prop_pre_post_inf_dist}, conditionally on $I_T=a$, the laws of the pre-$H_I$ and post-$H_I$ processes as well as  their supremums and  maximum drawdowns are characterized in terms of the scale functions. 
We identify the law of the pre-$H_I$ process under conditioning $H_I < H_S, I_T=a, S_T=b$ in Proposition \ref{prop:HI-before-HS}.
Finally in Theorem \ref{theorem:mdd-HIHS}, we find the conditional distribution of maximum drawdown of the pre-infimum process and recall for the intermediate and the post-supremum processes under the conditioning of Proposition \ref{prop:HI-before-HS}.

\section{Preliminaries}\label{sec:preliminaries}

Throughout, $X=(X_t)_{t\ge0}$ is a spectrally negative L\'{e}vy process with unbounded variation on a complete filtered probability space $(\Omega,\F,(\F_t)_{t\ge0},\Pp)$.  Its Laplace exponent $\psi$ is defined by
\begin{equation}\label{eq:laplace-exponent}
\E_x\big[e^{\lambda(X_t-X_0)}\big]=e^{t\psi(\lambda)},\qquad \lambda\ge0,
\nonumber\end{equation}
on the set where the expectation is finite where
$$
\Psi(\lambda)=\frac{\sigma^2}{2} \lambda^2+\mu \lambda+\int_{(0, \infty)}\left[e^{-\lambda y}-1+\lambda y\right] \Pi(\mathrm{d} y), \lambda \geq 0$$
 with a Lévy measure  $\Pi$ \text { of } $-X$ \text { satisfying }
$\int_{(0, \infty)}\left(y \wedge y^2\right) \Pi(\mathrm{d} y)<\infty$, drift parameter $\mu \in \mathbb{R}$ and diffusion parameter $\sigma\ge 0$.

%...$ yazalim ki $\sigma$ gozuksun.}
For $q\ge0$, let
\[
\Phi(q)=\sup\{\lambda\ge0:\psi(\lambda)=q\}.
\]
The $q$-scale function $W^{(q)}:\R\to[0,\infty)$ is characterized by $W^{(q)}(x)=0$ for $x<0$ and
\begin{equation}\label{eq:scale-def}
\int_0^\infty e^{-\lambda x}W^{(q)}(x)\,\dd x=\frac{1}{\psi(\lambda)-q},
\qquad \lambda>\Phi(q).\nonumber
\end{equation}
The associated second $q$-scale function is found as
\begin{equation}\label{eq:Z-def}
Z^{(q)}(x)=1+q\int_0^x W^{(q)}(y)\,\dd y,
\qquad x\ge0.\nonumber
\end{equation}
%Thus, 
%\begin{equation}\label{eq:Zprime}
%Z^{(q)\prime}(x)=qW^{(q)}(x).
%\end{equation}
Since  $X$ is assumed to have unbounded variation, the origin is regular for both half-lines and the scale function  $W^{(q)}$ is continuously differentiable on $(0,\infty)$ \cite[Lem.2.4]{Kuznetsov2013}. 

We consider the first passage times defined by
\[
\tau_c^+:=\inf\{t\ge0:X_t>c\},\qquad
\tau_c^-:=\inf\{t\ge0:X_t<c\}.
\]
For $a<x<b$, the standard two-sided exit identities are
\begin{eqnarray}
\E_x\left[e^{-q\tau_b^+};\tau_b^+<\tau_a^-\right]
&=&\frac{W^{(q)}(x-a)}{W^{(q)}(b-a)},\label{eq:two-sided-up}\\
\E_x\left[e^{-q\tau_a^-};\tau_a^-<\tau_b^+\right]
&=&Z^{(q)}(x-a)-Z^{(q)}(b-a)\frac{W^{(q)}(x-a)}{W^{(q)}(b-a)}.\label{eq:two-sided-down}
\end{eqnarray}
We also use the one-sided identity
\begin{equation}\label{eq:one-sided-down}
\E_x\left[e^{-q\tau_a^-};\tau_a^-<\infty\right]
=Z^{(q)}(x-a)-\frac{q}{\Phi(q)}W^{(q)}(x-a),\qquad x>a.\nonumber
\end{equation}
  In the following sections, the killed processes are formulated through passage times $\tau_c^\pm$ except for the pre-infimum process where the ordinary hitting time   $\tau_c:=\inf\{t\ge0:X_t=c\}$ is required. On the other hand, $\tau_c^{+} =\tau_c$ because $X$ has no positive jumps and the process creeps upwards.
  
Let $T\sim \operatorname{Exp}(\gamma)$, $\gamma>0$, be independent of $X$.  We define
\[
S_T:=\sup_{0\le t\le T}X_t,\qquad I_T:=\inf_{0\le t\le T}X_t,
\]
and denote by $H_S$ and $H_I$ the times at which these extrema are attained, respectively, that is, 
\[
H_S=\inf\{t < T: X_{t}=S_T\}, \qquad H_I=\inf\{t < T: X_t=I_T\} \; .
\]
 Since $X$ has no positive jumps, the  supremum is attained continuously, and also 0 is regular for $(0, \infty)$.
%{\it Theorem 3.1 of \cite{Millar1977}}: Let $T$ be a random time, let $I_T=\inf_{ 0\le t\le T} X_t $, and let $M$ be the point in $[0,T]$ where this minimum is achieved. Assume that $P\{M<T\}>0$. 
%Then, 
%\[
%i) P \{X_M = I_T|M<T\} = 1 \quad \mbox{if 0 is regular for $(0, \infty)$}
%\]
%\[
%ii) P \{X_M > I_T|M<T\} = 1 \quad \mbox{if 0 is not regular for $(0, \infty)$}
%\]
It follows from  \cite[Thm.3.1]{Millar1977} that 
\begin{equation}
X_{H_I} = I_T  \label{temp}
\end{equation}
almost surely under $H_I<T$. % where $H_I = \inf\{t<T: X_t=I_{T}\}$. 
    This can happen in two ways: the process $X$ either jumps into $I_T$ or $X$ is continuous at  $I_T$. If 0 is regular also for $(-\infty,0)$, which happens when $X$ is of unbounded variation, then $X$ is continuous at $I_T$ (\cite{Millar1977}, pg.370, Remark after Prop.2.4).  

%\begin{assumption}\label{ass:standing}
%All derivatives at the origin are interpreted as right derivatives.
    % Since it has a non-trivial Gaussian component $\sigma>0$, we have
    % \begin{equation}\label{eq:Wprime0}
    % W^{(q)\prime}(0+)=\frac{2}{\sigma^2}<\infty,\qquad q\ge0.
    % \end{equation}

%\end{assumption}

%\begin{remark}\label{rem:assumption}
%The Gaussian-component assumption is imposed to avoid separate entrance-law normalizations in cases where $W^{(q)\prime}(0+)=\infty$.  The same formulae can be rewritten under weaker hypotheses by replacing the finite constants $W^{(q)\prime}(0+)$ with the corresponding boundary limits.  We do not pursue this extension here.
%\end{remark}

%\begin{convention}[Conditional laws at  extrema]\label{conv:conditioning}
Statements involving conditioning on $I_T=a$ or on $(I_T,S_T)=(a,b)$ are interpreted as statements under fixed versions of the corresponding regular conditional laws. We use the notation
\begin{eqnarray}
\Pp_{a,b}^{<}(\,\boldsymbol\cdot\,)
:=\Pp_0(\,\boldsymbol\cdot\mid H_I<H_S,I_T=a,S_T=b) \label{convention}
\end{eqnarray}
for the conditional distribution with $a<0<b$.  
%The one-extremum conditioning $\Pp_0(\,\boldsymbol\cdot\mid I_T=a)$ is understood analogously.  The formulae below identify these versions through density-ratio limits.
%\end{convention}
A positive excessive function $h$ for a killed Markov process defines the Doob $h$-transform by
\begin{equation*}\label{eq:h-transform}
P_t^h(x,\dd y)=\frac{h(y)}{h(x)}\Pp_x\{X_t\in\dd y,\ t<\zeta\},
\end{equation*}
where $\zeta$ is the killing time.  Positivity and excessivity of the functions $h$  follow from their representation as conditional  entrance probabilities for the corresponding killed semigroups.

\section{Path decomposition at the  infimum}\label{sec:infimum-decomposition}

The following proposition identifies the law of the pre-infimum and post-infimum process and the distribution of the supremum within these parts of the path given the value of the infimum.
\begin{theorem}\label{Prop_pre_post_inf_dist}
	 Conditionally on $I_{T}=a$, the pre-$H_{I}$ and post-$H_{I}$ processes are independent. 
     
i. The law of the pre-$H_I$ process is an $h$-transform of the law of the spectrally negative L\'{e}vy process killed at  $\tau_a\wedge T $ with 
\[
h(x) = \frac{\gamma}{\Phi(\gamma)} W^{(\gamma)\prime}(x-a) -\gamma W^{(\gamma)}(x-a)
\] 
for $x\ge a$. Moreover, the distribution of the supremum and the maximum drawdown of the pre-infimum process, $S_{H_{I}}$, $M_{0, H_{I}}^{-},$ respectively, when $I_{T}=a$ are given by
	\begin{eqnarray*}
		% &&\mathbb{P}_0\left(S_{H_{I}} < b \mid I_{T}=a\right)=\frac{ \frac{\gamma}{\Phi(\gamma)} W^{(\gamma)\prime}(0) }{\frac{\gamma}{\Phi(\gamma)} W^{(\gamma)\prime}(-a) -\gamma W^{(\gamma)}(-a)}\left[Z^{(\gamma)}(-a) - Z^{(\gamma)}(b-a) \frac{W^{(\gamma)}(-a)}{W^{(\gamma)}(b-a)}\right]
        &&\mathbb{P}_0\left\{S_{H_{I}} < b \mid I_{T}=a\right\}=1-\frac{ \frac{\gamma}{\Phi(\gamma)} W^{(\gamma)\prime}(b-a) -\gamma W^{(\gamma)}(b-a) }{\frac{\gamma}{\Phi(\gamma)} W^{(\gamma)\prime}(-a) -\gamma W^{(\gamma)}(-a)}\left[\frac{W^{(\gamma)}(-a)}{W^{(\gamma)}(b-a)}\right]
	\end{eqnarray*}
	for  $b\ge 0$, and
    \begin{eqnarray*}
\Pp\left\{M_{0, H_{I}}^{-}< d \mid I_{T}=a \right\} =1-\frac{ \frac{\gamma}{\Phi(\gamma)} W^{(\gamma)\prime}(d) -\gamma W^{(\gamma)}(d) }{\frac{\gamma}{\Phi(\gamma)} W^{(\gamma)\prime}(-a) -\gamma W^{(\gamma)}(-a)}\left[\frac{W^{(\gamma)}(-a)}{W^{(\gamma)}(d)}\right]
\end{eqnarray*}
for $d\ge -a.$

ii. The law of the post-$H_{I}$ process is the $h$-transform of the law of the spectrally negative L\'{e}vy process killed at $T \wedge \tau_{a}^{-}$with	
	$$
	h(x)=1-Z^{(\gamma)}(x-a)+\frac{\gamma}{\Phi(\gamma)} W^{(\gamma)}(x-a).
	$$	
	Moreover, the distribution of the supremum of the post-infimum process, $S_{H_{I}, T}$, when $I_{T}=a$ is given by
%\begin{equation}
	%\mathbb{P}_{a}\left(S_{H_{I}, T} \geq b \mid I_{T}=a\right)=\frac{1-Z^{(\gamma)}(b-a)+\frac{\gamma}{\Phi(\gamma)} W^{(\gamma)}(b-a)}{\frac{\gamma}{\Phi(\gamma)}W^{(\gamma)}(b-a)}
%\end{equation}
\[
\mathbb{P}_{a}\left\{S_{H_{I}, T} \le b \mid I_{T}=a\right\}  = \frac{\Phi(\gamma)({Z^{(\gamma)}}(b-a)-1)}{\gamma W^{(\gamma)}(b-a)} 
\]
	for $b>a$.
\end{theorem}

Proof: 
i. Recall $\tau_a = \inf\{t\ge 0: X_t=a\}$. Note that $0<\tau_a<T$ when 
0 is regular  for both $(0,\infty)$ and  $(-\infty,0)$ by \cite[Prop.2.1]{Millar1977}. Given $\{I_T =a\}$, we have $X_{H_I}=a$ by \eqref{temp} and 
also because  the  minimum before $T$ occurs at a single point in time almost surely, that is, $H_I=\tau_a$, by \cite[Prop.2.2]{Millar1977}. These justify the use of $\tau_a$ and $H_I$ interchangeably below.  For $a<0$ and $x,y>a$, we have
\begin{eqnarray*}
\mathbb{P}_x \left\{X_t \in {d} y, t<\tau_a \mid I_T=a\right\}&=&\frac{ \frac{ d}{ d a}\mathbb{P}_x \left\{X_t \in d y, t<H_I , I_T\leq a\right\}}{\frac{d}{ d a} \mathbb{P}_x\left\{I_T \leqslant a \right\}}\\
&=&\frac{ {\frac{d}{da}}\mathbb{P}_x \left\{I_T \leqslant a \mid X_t \in d y, t<\tau_a\wedge T  \right\}\mathbb{P}_x \left\{X_t \in d y, t<\tau_a {\wedge T}  \right\}}{\frac{d}{ d a} \mathbb{P}_x\left\{I_T \leqslant a \right\}}\\ 
&=&\frac{\frac{ d}{d a} \mathbb{P}_y\left\{I_T \leqslant a \right\}\mathbb{P}_x \left\{X_t \in d y, t<\tau_a { \wedge T}   \right\}}{\frac{ d}{ d a} \mathbb{P}_x\left\{I_T \leqslant a \right\}}\\ 
&=:&\frac{ h(y) } { h(x) } \:\mathbb{P}_x \left\{X_t \in {d} y, t< \tau_a { \wedge T} \right\}  
\end{eqnarray*}
where under $\mathbb{P}_x$, for  $x>a$, we have
\begin{eqnarray*}
h(x) 
%&=& \frac{\mathrm{d}}{\mathrm{d}a} \mathbb{P}_x\{I_T\le a\} = \frac{\mathrm{d}}{\mathrm{d}a} \mathbb{P}_{-x}\{-I_T\le -a\} =
% \frac{\mathrm{d}}{\mathrm{d}a} \mathbb{P}\{-I_T\le x-a\} 
&=& \frac{{d}}{{d}a} \mathbb{P}_x\{I_T\le a\} =\frac{{d}}{{d}a} \mathbb{P}\{I_T\le a-x\} = \frac{{d}}{{d}a} \mathbb{P} \{-I_T> x-a\}  \\
&=& \frac{{d}}{{d}a} \int_{x-a}^{\infty} \left[ \frac{\gamma}{\Phi(\gamma)} W^{(\gamma)}( dy) -\gamma W^{(\gamma)}(y) \: {dy} \right]  \\
&=& \frac{{d}}{{d}a} \int_{x-a}^{\infty} \left[ \frac{\gamma}{\Phi(\gamma)} W^{(\gamma)\prime}(y)  -\gamma W^{(\gamma)}(y) \right]\:{dy}  \\
&=& \frac{\gamma}{\Phi(\gamma)} W^{(\gamma)\prime}(x-a) -\gamma W^{(\gamma)}(x-a) 
\end{eqnarray*} 
by \cite[Eq.8.24]{Kyprianou2014}  since we have assumed that L\'{e}vy process $X$ is of unbounded variation so that $W^{(\gamma)}$  is differentiable by  \cite[Lem.2.4]{Kuznetsov2013}.
Therefore, the law of the pre-$H_I$ process is an $h$-transform of the law of the spectrally negative L\'{e}vy process killed at  $\tau_a\wedge T $ with 
\[
h(x) = \frac{\gamma}{\Phi(\gamma)} W^{(\gamma)\prime}(x-a) -\gamma W^{(\gamma)}(x-a).
\] 
Then, it follows that
\begin{eqnarray*}\label{eqnpreinfsupremumgreaterthanb}
\mathbb{P}_0\left\{S_{H_{I}} < b \mid I_{T}=a\right\} &=& 1-\mathbb{P}_0\left\{S_{H_{I}} \geq b \mid I_{T}=a\right\}\nonumber\\&=&1-\dfrac{h(b)}{h(0)}\mathbb{P}_0\left(\tau_b^+<\tau_a^-<, \tau_b^+<T\right)\nonumber\\ &=& 1-\dfrac{h(b)}{h(0)}\left[ \frac{W^{(\gamma)}(-a)}{W^{(\gamma)}(b-a)}\right]\nonumber\\ &
=&1-\frac{ \frac{\gamma}{\Phi(\gamma)} W^{(\gamma)\prime}(b-a) -\gamma W^{(\gamma)}(b-a) }{\frac{\gamma}{\Phi(\gamma)} W^{(\gamma)\prime}(-a) -\gamma W^{(\gamma)}(-a)}\left[\frac{W^{(\gamma)}(-a)}{W^{(\gamma)}(b-a)}\right]
\end{eqnarray*}
This result satisfies the correct boundary limits at $b=0$ and as $b\rightarrow \infty$ in view of $\lim_{x\rightarrow \infty} W^{(\gamma)\prime}(x)/W^{(\gamma)}(x) =\Phi(\gamma).$ Then, we also get the distribution of $M_{0, H_{I}}^{-}$ by the identity
% The complementary probability is given by
% \begin{eqnarray}\label{complement}
% \mathbb{P}_0\left(S_{H_{I}} > b \mid I_{T} = a\right) &=& \lim _{x \rightarrow 0}\dfrac{h(b)}{h(x)}\mathbb{P}_x\left(\tau_a^->\tau_b^+, \tau_b^+<T\right)\nonumber\\ &= & \lim _{x \rightarrow 0}\dfrac{h(b)}{h(x)}     \frac{W^{(\gamma)}(x-a)}{W^{(\gamma)}(b-a)} \nonumber \\ \nonumber\\ &
% = &\lim _{x \rightarrow 0}\frac{ \frac{\gamma}{\Phi(\gamma)} W^{(\gamma)\prime}(b-a) -\gamma W^{(\gamma)}(b-a)}{\frac{\gamma}{\Phi(\gamma)} W^{(\gamma)\prime}(x-a) -\gamma W^{(\gamma)}(x-a)} \: \frac{W^{(\gamma)}(x-a)}{W^{(\gamma)}(b-a)}  \\ \nonumber \\
% &=&\frac{ \frac{\gamma}{\Phi(\gamma)} W^{(\gamma)\prime}(b-a) -\gamma W^{(\gamma)}(b-a)}{\frac{\gamma}{\Phi(\gamma)} W^{(\gamma)\prime}(-a) -\gamma W^{(\gamma)}(-a)} \: \frac{W^{(\gamma)}(-a)}{W^{(\gamma)}(b-a)} 
% \end{eqnarray}
\begin{eqnarray*}
\Pp\left\{M_{0, H_{I}}^{-}< d \mid I_{T}=a \right\} = \mathbb{P}_0\left\{S_{H_{I}} < d+a \mid I_{T}=a\right\}.
\end{eqnarray*}

ii.  Proof was given in \cite{VardarAcar2021}.

\begin{remark}(Comparison of the pre-$H_I$ process with \cite[Thm.3.2,Cor.3.4]{Salminen2007}) 
In Proposition~\ref{Prop_pre_post_inf_dist}, conditioned on $I_T=a$, the pre-$H_I$ process 
is the Doob-$h$ transform of the spectrally negative L\'{e}vy process killed at $\tau_a \wedge T$, 
with   
\[
h(x)=\frac{\gamma}{\Phi(\gamma)} W^{(\gamma)\prime}(x-a)-\gamma W^{(\gamma)}(x-a).
\]
For  Brownian motion with no drift, the scale function is given by 
$W^{(\gamma)}(x)=\tfrac{2}{\sqrt{2\gamma}}\sinh(\sqrt{2\gamma}x)$ and its derivative is
$W^{(\gamma)\prime}(x)=2\cosh(\sqrt{2\gamma}x)$. Substituting these into $h$ gives 
\[
h(x)=\sqrt{2\gamma}\big(\cosh(\sqrt{2\gamma}(x-a))-\sinh(\sqrt{2\gamma}(x-a))\big)
   =\sqrt{2\gamma}\, e^{-\sqrt{2\gamma}(x-a)}.
\]
The Doob-$h$ transform with $h(x)=\sqrt{2\gamma}\,e^{-\sqrt{2\gamma}(x-a)}$\ changes the generator of the Brownian 
motion from $\mathcal{L}f=\tfrac12 f''$ to 
\[
\mathcal{L}^h f
 =\frac{1}{h}\mathcal{L}(hf)
 =\tfrac12 f''-\sqrt{2\gamma}f',
\]
which is the generator of a Brownian motion with drift $-\sqrt{2\gamma}$ killed upon hitting $a$. 
This coincides exactly with the result of Salminen \& Vallois (\cite{Salminen2007}, Theorem 3.2). 
Also in the Brownian case, we have
\begin{eqnarray}\mathbb{P}_0\left(S_{H_{I}} < b \mid I_{T}=a\right)
&=&1+\frac{e^{-b\sqrt{2\gamma}}\sinh{(a\sqrt{2\gamma})}}{\sinh{((b-a)\sqrt{2\gamma})}} \nonumber\\
&=&\frac{e^{-a\sqrt{2\gamma}}\sinh{(b\sqrt{2\gamma})}}{\sinh{((b-a)\sqrt{2\gamma})}} \label{Br}
\end{eqnarray} which confirms the result given  in \cite[ Cor.3.4]{Salminen2007}.  
\end{remark}

\section{Path decomposition under $H_I<H_S$}\label{sec:ordered-decomposition}

We next condition on both  extrema and on the ordering that the infimum occurs before the supremum.  We fix $a<0<b$, and use the notation $\Pp_{a,b}^{<}$ of \eqref{convention}.

\begin{proposition}[Decomposition under $H_I<H_S$]\label{prop:HI-before-HS}
Under $\Pp_{a,b}^{<}$, the pre-$H_I$ process, the intermediate process
$$
\{X_{H_I+u}:0\le u\le H_S-H_I\},
$$
and the post-$H_S$ process are independent.  Their laws are as follows.
\begin{enumerate}[label=\textup{(\roman*)}]
\item The pre-infimum process is the Doob transform of the process killed on leaving $(a,b)$ through the lower boundary, with
\begin{equation}\hspace{-1,4cm}\label{eq:h1-HIHS}
h_1(x)=
\frac{W^{(\gamma)\prime}(x-a)W^{(\gamma)}(b-a)-W^{(\gamma)\prime}(b-a)W^{(\gamma)}(x-a)}
     {\big(W^{(\gamma)}(b-a)\big)^2},
\qquad a<x<b.
\end{equation}
%The lower boundary value is understood as
%\begin{equation}\label{eq:h1-lower}
%h_1(a)=\frac{W^{(\gamma)\prime}(0+)}{W^{(\gamma)}(b-a)} = \frac{2/\sigma^2 }{\, W^{(\gamma)}(b-a)}.
%\end{equation}

\item The intermediate process has the law of $a+\Gamma$, where $\Gamma$ is the Doob transform of a spectrally negative L\'{e}vy process killed at $\tau_{b-a}^+\wedge\tau_0^-$, with
\begin{equation}\label{eq:h2-HIHS}
h_2(z)=e^{-\Phi(\gamma)(b-a-z)}\frac{W^{(\gamma)}(z)}{W^{(\gamma)}(b-a)},
\qquad 0<z<b-a.
\end{equation}
Equivalently, in the original coordinate $x=a+z$,
\[
h_2^a(x)=e^{-\Phi(\gamma)(b-x)}\frac{W^{(\gamma)}(x-a)}{W^{(\gamma)}(b-a)}.
\]

\item The post-$H_S$ process is the Doob transform of the process killed on leaving $(a,b)$ before $T$, with
\begin{equation}\label{eq:h3-HIHS}
h_3(x)=1-Z^{(\gamma)}(x-a)+\big(Z^{(\gamma)}(b-a)-1\big)\frac{W^{(\gamma)}(x-a)}{W^{(\gamma)}(b-a)},
\qquad a<x<b.
\end{equation}
\end{enumerate}
\end{proposition}

%\begin{proof}

Proof: By \cite[Lem.1]{VardarAcar2021}  the pre-$H_{S}$ process is a L\'{e}vy process with Laplace exponent $\widetilde{\psi}$ given by $
\widetilde\psi(\lambda)=\psi(\lambda+\Phi(\gamma))-\gamma$.  It is governed by Esscher-transformed law with parameter $\Phi(\gamma)$.
 Let $\Pp^{\Phi(\gamma)}$ denote the probability corresponding to this change of measure.

The infimum of the presupremum process is also equal to $a$. Under the given conditions, note that the intermediate process  is independent from the  preinfimum process  and  both governed by  $\Pp^{\Phi(\gamma)}$.  
Therefore, the preinfimum process starts at $0$ and evolves independently from the other parts  as a process killed at the first hitting time at $a$ before the exponential time $T$ under the condition that it stays below $b$, and hits $a$ before going above $b$. 

For $  a<x,y<b$, we have
\begin{eqnarray*}  
P_t(x,dy) & = &  \Pp_x \{X_t\in dy, t<H_I\, |\, H_I<H_S, I_T =a, S_T=b\} \\
& = & \frac{\frac{d}{da} \Pp_x \{X_t\in dy, t<H_I, H_I<H_S, I_T \leq a \, | \,S_T=b\} }{\frac{d}{da}\Pp_x \{ H_I<H_S, I_T \leq a \, | \,S_T=b\}     }  \\
& = & \frac{ \frac{d}{da}\Pp_y^{\Phi(\gamma)} \{ I_{\tau_b^+} \leq a \, | X_t\in dy, t<\tau_a
\} \Pp_x\{X_t\in dy, t<\tau_a\}}{\frac{d}{da}\Pp_x^{\Phi(\gamma)} \{  I_{\tau_b^+} \leq a \} }\\
& = & \frac{\frac{d}{da} \Pp_y^{\Phi(\gamma)} \{ I_{\tau_b^+} \leq a  \} }{\frac{d}{da}\Pp_x^{\Phi(\gamma)} \{ I_{\tau_b^+} \leq a \}     } \Pp_x^{\Phi(\gamma)} \{X_t\in dy, t<\tau_a\}\\
&=&\frac{\frac{d}{da} \Pp_y^{\Phi(\gamma)} \{ I_{\tau_b^+} \leq a  \} }{\frac{d}{da}\Pp_x^{\Phi(\gamma)} \{ I_{\tau_b^+} \leq a \} } \frac{e^{-\Phi(\gamma)(b-y)}}{e^{-\Phi(\gamma)(b-x)}} \mathbb{P}_{x}\left\{X_{t} \in {d} y, t<\tau_{a}\right\} \\
&=:& \frac{h_1(y)}{h_1(x)} \,\mathbb{P}_{x}\left\{X_{t} \in {d} y, t<\tau_{a} \right\} 
\end{eqnarray*}
where
$$\frac{d}{da}\Pp_x^{\Phi(\gamma)} \{ I_{\tau_b^+} \leq a \}=\frac{d}{da}\Pp_x \{ H_I<H_S, I_{T} \leq a \, | \,S_T=b\} .$$

The scale function under $\Pp_x^{\Phi(\gamma)}$ is
\[
W_{\Phi(\gamma)}(u)=e^{-\Phi(\gamma)u}W^{(\gamma)}(u)
\]
for $u\ge 0$. Since 
\[
\Pp_x^{\Phi(\gamma)}\{I_{\tau_b^+}\le a\}
=1-\frac{W_{\Phi(\gamma)}(x-a)}{W_{\Phi(\gamma)}(b-a)},
\]
%By \cite{VardarAcar2017}, the result given in the proof [Theorem 1, ii.)] we have
it follows that
$$
1-\mathbb{P}_{x}^{\Phi(\gamma)}\left\{I_{\tau_b^{+}} \leq a\right\}=e^{-\Phi(\gamma)(x-b)}\frac{W^{(\gamma)}(x-a)}{W^{(\gamma)}(b-a)}.
$$ 
and
%$$
%\mathbb{P}^{\Phi(\gamma)}_x\left(I_{\tau_b^{+}} \leq a\right) = \mathbb{P}^{\Phi(\gamma)}_0\left(I_{\tau_{b-x}^{+}} \leq a-x\right)=1-e^{-\phi(\gamma)(x-b)}\frac{W^{(\gamma)}(x-a)}{W^{(\gamma)}(b-a)}.
%$$ 
\begin{eqnarray}&&\hspace*{-1,5cm} \frac{d}{da}\mathbb{P}^{\Phi(\gamma)}_x\left\{I_{\tau_b^{+}} \leq a\right\}=\frac{W^{\Phi(\gamma)\prime}(x-a) W^{\Phi(\gamma)}(b-a)-W^{\Phi(\gamma)\prime}(b-a) W^{\Phi(\gamma)}(x-a)}{(W^{\Phi(\gamma)}(b-a))^{2}}\nonumber\\
&=&\frac{[\Phi(\gamma)e^{-\Phi(\gamma)(x-a)}W^{(\gamma)}(x-a)+W^{(\gamma)\prime}(x-a)e^{-\Phi(\gamma)(x-a)}] e^{-\Phi(\gamma)(b-a)}W^{(\gamma)}(b-a)}{(W^{\Phi(\gamma)}(b-a))^{2}}\nonumber\\
&& -\frac{[\Phi(\gamma)e^{-\Phi(\gamma)(b-a)}W^{(\gamma)}(b-a)+W^{(\gamma)\prime}(b-a)e^{-\Phi(\gamma)(b-a)}] e^{-\Phi(\gamma)(x-a)}W^{(\gamma)}(x-a)}{(W^{\Phi(\gamma)}(b-a))^{2}}\nonumber\\
&=&\frac{e^{-\Phi(\gamma)(x-a)}W^{(\gamma)}(b-a)W^{(\gamma)\prime}(x-a)-e^{-\Phi(\gamma)(x-a)}W^{(\gamma)}(x-a)W^{(\gamma)\prime}(b-a)}{e^{-\Phi(\gamma)(b-a)}(W^{(\gamma)}(b-a))^{2}}\nonumber\\
&=&e^{-\Phi(\gamma)(x-b)}\frac{W^{(\gamma)}(b-a)W^{(\gamma)\prime}(x-a)-W^{(\gamma)}(x-a)W^{(\gamma)\prime}(b-a)}{(W^{(\gamma)}(b-a))^{2}}
\end{eqnarray}
Thus, we get
$$h_1(x)=\frac{W^{(\gamma)}(b-a)W^{(\gamma)\prime}(x-a)-W^{(\gamma)}(x-a)W^{(\gamma)\prime}(b-a)}{(W^{(\gamma)}(b-a))^{2}}$$

    %Condition first on $S_T=b$.  

Conditional independence of the three pieces follows by applying the Markov property at the two extremes \cite{Millar1977}.

ii. and iii. Proofs were given in \cite{VardarAcar2021}.
%\end{proof}

\begin{remark}
    Note that taking the derivative with respect to $a$ and returning to the original coordinates gives the boundary derivative
\[
\frac{d}{da}\Pp_x^{\Phi(\gamma)}\{I_{\tau_b^+}\le a\}
=e^{-\Phi(\gamma)(x-b)}h_1(x).
\]
The Esscher factors cancel in the density-ratio kernel, yielding the $h_1$-transform in \eqref{eq:h1-HIHS}.  Since it is the derivative of a hitting probability with respect to the lower boundary, $h_1$ is positive and excessive for the killed semigroup.  
For the intermediate component, after the path has entered at $a$ it is conditioned to reach $b$ before returning below $a$ and before exponential killing.  The excessive function is therefore the discounted two-sided upward exit probability under the pre-supremum Esscher transform.  Applying \eqref{eq:two-sided-up} gives $h_2$-transform in \eqref{eq:h2-HIHS}.
\end{remark}

The next theorem gives the maximum-drawdown laws of the three independent pieces stated in Proposition~\ref{prop:HI-before-HS} for an unbounded variation L\'{e}vy process.  The cases outside the non-trivial range are stated explicitly. 

\begin{theorem}[Maximum drawdowns under $H_I<H_S$]\label{theorem:mdd-HIHS}
Fix $a<0<b$ and let $r=b-a$.  Under $\Pp_{a,b}^{<}$, the following identities hold.
\begin{enumerate}[label=\textup{(\roman*)}]
\item For the pre-infimum process,
\begin{equation*}\hspace{-1cm}\small\label{eq:f1-piecewise}
\Pp_{a,b}^{<}\{M_{0,H_I}^-<d\}=\begin{cases}
0, & d\le -a,\\
\displaystyle
1-\left[\frac{W^{(\gamma)\prime}(d)W^{(\gamma)}(r)-W^{(\gamma)\prime}(r)W^{(\gamma)}(d)}{W^{(\gamma)\prime}(-a)W^{(\gamma)}(r)-W^{(\gamma)\prime}(r)W^{(\gamma)}(-a) }\right]
\left(\frac{W^{(\gamma)}(-a)}{W^{(\gamma)}(d)}\right),
& -a<d<r,\\
1, & d\ge r.
\end{cases}
\end{equation*}

\item For the intermediate process,
\begin{equation*}\label{eq:f2-piecewise}
\Pp_{a,b}^{<}\{M_{H_I,H_S}^-<d\}=\begin{cases}
0, & d\le0,\\[4pt]
\displaystyle
\frac{W^{(\gamma)}(r)}{W^{(\gamma)}(d)}
\exp\left\{-(r-d)\frac{W_+^{(\gamma)\prime}(d)}{W^{(\gamma)}(d)}\right\},
& 0<d<r,\\[12pt]
1, & d\ge r.
\end{cases}
\end{equation*}

\item For the post-supremum process,
\begin{equation*}\label{eq:f3-piecewise}
\Pp_{a,b}^{<}\{M_{H_S,T}^-<d\}=\begin{cases}
0, & d\le0,\\[4pt]
\displaystyle
\frac{\big(Z^{(\gamma)}(d)-1\big)\dfrac{W^{(\gamma)\prime}(d)}{W^{(\gamma)}(d)}-\gamma W^{(\gamma)}(d)}
     {\big(Z^{(\gamma)}(r)-1\big)\dfrac{W^{(\gamma)\prime}(r)}{W^{(\gamma)}(r)}-\gamma W^{(\gamma)}(r)},
& 0<d<r,\\[14pt]
1, & d\ge r.
\end{cases}
\end{equation*}
\end{enumerate}
\end{theorem}

%\begin{proof}

Proof: (i) Under $\Pp_{a,b}^{<}$ the pre-infimum path is governed by the $h_1$-transform in \eqref{eq:h1-HIHS}.  Because the path starts at $0$ and ends at the  infimum level $a$, the event $\{M_{0,H_I}^-<d\}$ is impossible when $d\le -a$.  If $-a<d<r$, no drawdown of size $d$ can occur before the path reaches $a$ precisely when the path exits the interval $(a,a+d)$ through the lower boundary before crossing above $a+d$.  Hence, we have
\begin{eqnarray*}
\Pp_{a,b}^{<}\{M_{0,H_I}^-<d\}&=&1-\Pp\left\{M_{0, H_{I}}^{-}\geq d \mid H_{I}<H_{S}, I_{T}=a, S_{T}=b\right\}=:1-\Pp_0^{h_{1}}(\tau_{a+d}^{+}<\tau_a^- )\nonumber\\ &=& 1-\frac{h_1(a+d)}{h_1(0)}\Pp_0(\tau_{a+d}^{+}<\tau_a^-, \tau_{a+d}^{+}<T)\nonumber\\
&=&1-\frac{h_1(a+d)}{h_1(0)}\left(\frac{W^{(\gamma)}(-a)}{W^{(\gamma)}(d)}\right)
\end{eqnarray*}
where $h_1$ is given in \eqref{eq:h1-HIHS}.
%\begin{eqnarray*}
%    h(a)&=&\left[\frac{W^{(\gamma)\prime}(0)W^{(\gamma)}(b-a)-W^{(\gamma)\prime}(b-a)W^{(\gamma)}(0)}{W^{(\gamma)}(b-a)^2}\right]
%\end{eqnarray*}
%\begin{eqnarray*}
%    h(0)&=&\left[\frac{W^{(\gamma)\prime}(-a)W^{(\gamma)}(b-a)-W^{(\gamma)\prime}(b-a)W^{(\gamma)}(-a)}{W^{(\gamma)}(b-a)^2}\right]
%\end{eqnarray*}
% and
%  \begin{eqnarray*}
%     \frac{h_1(a)}{h_1(0)}&=&\left[\frac{W^{(\gamma)\prime}(0)W^{(\gamma)}(b-a)}{W^{(\gamma)\prime}(-a)W^{(\gamma)}(b-a)-W^{(\gamma)\prime}(b-a)W^{(\gamma)}(-a) }\right]
% \end{eqnarray*} 
% {\color{red} Alternatively,}
% \begin{eqnarray*}
% \Pp_{a,b}^{<}\{M_{0,H_I}^->d\}&=&\Pp\left\{M_{0, H_{I}}^{-}>d \mid H_{I}<H_{S}, I_{T}=a, S_{T}=b\right\}=\Pp_0^{h_{1}}(\tau_a^- > \tau_{a+d}^{+})\nonumber\\ &=& \frac{h_1(a+d)}{h_1(0)}\Pp_0(\tau_a^- > \tau_{a+d}^{+},\tau_{a+d}^+<T)\nonumber\\
% &=&\frac{h_1(a+d)}{h_1(0)} \, \frac{W^{(\gamma)}(-a)}{W^{(\gamma)}(d)} 
% \end{eqnarray*}
 Therefore, we get
\begin{eqnarray*}
&&\hspace{-1cm}\Pp_{a,b}^{<}\{M_{0,H_I}^-<d\}=1-\Pp\left\{M_{0, H_{I}}^{-}\geq d \mid H_{I}<H_{S}, I_{T}=a, S_{T}=b\right\} \nonumber\\&=&1-\left[\frac{W^{(\gamma)\prime}(d)W^{(\gamma)}(b-a)-W^{(\gamma)\prime}(b-a)W^{(\gamma)}(d)}{W^{(\gamma)\prime}(-a)W^{(\gamma)}(b-a)-W^{(\gamma)\prime}(b-a)W^{(\gamma)}(-a) }\right]
\left(\frac{W^{(\gamma)}(-a)}{W^{(\gamma)}(d)}\right)\, .
\end{eqnarray*}
%Note that when $d=-a$, the probability is zero as the infimum equals $a$. 
%{\color{red} Buraya $h_1(a+)$ 'nin tam a'da kill etme islevi oldugu yazilacak.}

%???????????????????
%We prove the three formulae in the non-trivial ranges; the boundary cases follow from the deterministic ranges of the corresponding path pieces.
%???????????????????
%FAZLA KISA KANIT :For (i), under $\Pp_{a,b}^{<}$ the pre-infimum path is governed by the $h_1$-transform in \eqref{eq:h1-HIHS}.  Because the path starts at $0$ and ends at the  infimum level $a$, the event $\{M_{0,H_I}^-<d\}$ is impossible when $d\le -a$.  If $-a<d<r$, no drawdown of size $d$ can occur before the path reaches $a$ precisely when the path exits the interval $(a,a+d)$ through the lower boundary before crossing above $a+d$.  Hence
%\begin{align*}
%\Pp_{a,b}^{<}\{M_{0,H_I}^-<d\}
%&=\frac{h_1(a+)}{h_1(0)}
%\Pp_0\{\tau_a^-<\tau_{a+d}^+,\ \tau_a^-<T\}.
%\end{align*}
%Using \eqref{eq:two-sided-down} on the interval $(a,a+d)$ gives
%\[
%\Pp_0\{\tau_a^-<\tau_{a+d}^+,\ \tau_a^-<T\}
%=Z^{(\gamma)}(-a)-Z^{(\gamma)}(d)\frac{W^{(\gamma)}(-a)}{W^{(\gamma)}(d)}.
%\]
%The identity \eqref{eq:h1-lower} and the value of $h_1(0)$ yield the prefactor in \eqref{eq:f1-piecewise}.  If $d\ge r$, the pre-infimum drawdown is bounded above by the  range $r$, hence the probability is one by continuity of the strict-threshold formula.

(ii) Let
\[
\alpha_d:=\inf\{t\ge0:S_t-X_t>d\}
\]
be the first time at which the drawdown exceeds $d$.
 For a generic spectrally negative L\'{e}vy process, we have
\begin{eqnarray}\label{eqn13}
\Pp_{x}\left\{\tau_{b}^{+}<\alpha_{d}, \tau_{b}^{+}<\tau_{a}^{-}, \tau_{b}^{+}<T\right\} & =&\mathbb{E}_{x}\left[e^{-q \tau_{b}^{+}} ; \tau_{b}^{+} \leq \alpha_{d} \wedge \tau_{a}^{-}\right] \\
& = &\frac{W^{(\gamma)}(x-a)}{W^{(\gamma)}(d)} e^{-(b-d-a) \frac{W^{(\gamma)\prime}(d)}{W^{(\gamma)}(d)}} \nonumber 
\end{eqnarray}
by identity (\cite{Avram2019}, Eq.20), top row and middle column, given in terms of the drawdown time $\alpha_{d}$. This is the usual reflected-process exit identity; see, for example, \cite{Mijatovic2012} and \cite{VardarAcar2017}. In \cite[Thm.2i]{VardarAcar2021}, the  probability of maximum drawdown for the intermediate process was found in terms of the distribution of the pre-supremum process but contains an error. Therefore, the correct expression with $h_2$  is derived as
\begin{eqnarray}\hspace{-1cm}\small
 \Pp\left\{M_{H_{I}, H_{S}}^{-}<d \mid H_{I}<H_{S}, I_{T}=a, S_{T}=b\right\} &=&\lim _{x \rightarrow a} \frac{ \Pp_{x}^{h_2}\left\{\tau_{b}^{+}<\alpha_{d}, \tau_{b}^{+}<\tau_{a}^{-}, \tau_{b}^{+}<T\right\}}{\Pp_{x}^{h_2}\left\{\tau_{b}^{+}<\tau_{a}^{-} \wedge T\right\}} \nonumber\\
    &=&\lim _{x \rightarrow a} \frac{\frac{h_2(b)}{h_2(x)}\Pp_{x}\left\{\tau_{b}^{+}<\alpha_{d}, \tau_{b}^{+}<\tau_{a}^{-}, \tau_{b}^{+}<T\right\}}{\frac{h_2(b)}{h_2(x)}\Pp_{x}\left\{\tau_{b}^{+}<\tau_{a}^{-} \wedge T\right\}} \nonumber\\
	&=&\lim _{x \rightarrow a}\frac{W^{(\gamma)}(b-a)}{W^{(\gamma)}(x-a)}\frac{W^{(\gamma)}(x-a)}{W^{(\gamma)}(d)} e^{-(b-a-d)\left(\frac{W^{(\gamma)\prime}(d)}{W^{(\gamma)}(d)}\right)}\nonumber\\
&=&\frac{W^{(\gamma)}(b-a)}{W^{(\gamma)}(d)} e^{-(b-a-d)\left(\frac{W^{(\gamma)\prime}(d)}{W^{(\gamma)}(d)}\right)}
\end{eqnarray}
for $0<d<b-a.$ 
%The standard drawdown first-passage identity for spectrally negative L\'{e}vy processes gives, for $a<x<b$ and $0<d<r$,
%\begin{equation}\label{eq:drawdown-first-passage}
%\Pp_x\{\tau_b^+<\alpha_d,\ \tau_b^+<\tau_a^-,\ \tau_b^+<T\}
%=\frac{W^{(\gamma)}(x-a)}{W^{(\gamma)}(d)}
%\exp\left\{-(b-a-d)\frac{W_+^{(\gamma)\prime}(d)}{W^{(\gamma)}(d)}\right\}.
%\end{equation}
The boundary values are immediate from the fact that the intermediate path stays in $[a,b]$ and has range $r$.
%Proof was given in \cite{VardarAcar2021}.
%For (iii), under $\Pp_{a,b}^{<}$ the post-supremum path is governed by the $h_3$-transform.  For $0<d<r$,
%\begin{align*}
%\Pp_{a,b}^{<}\{M_{H_S,T}^-<d\}
%&=\lim_{x\uparrow b}
%\frac{\Pp_x\{T<\tau_{b-d}^-\wedge\tau_b^+\wedge\tau_a^-\}}
%     {\Pp_x\{T<\tau_b^+\wedge\tau_a^-\}}.
%\end{align*}
%Since $a<b-d$, the numerator equals $\Pp_x\{T<\tau_{b-d}^-\wedge\tau_b^+\}$.  Applying \eqref{eq:two-sided-up}--\eqref{eq:two-sided-down} to the numerator and denominator and then using L'Hospital's rule as $x\uparrow b$ gives
%\[
%\frac{\big(Z^{(\gamma)}(d)-1\big)\dfrac{W^{(\gamma)\prime}(d)}{W^{(\gamma)}(d)}-Z^{(\gamma)\prime}(d)}
%     {\big(Z^{(\gamma)}(r)-1\big)\dfrac{W^{(\gamma)\prime}(r)}{W^{(\gamma)}(r)}-Z^{(\gamma)\prime}(r)}.
%\]
    % Finally, \eqref{eq:Zprime} gives $Z^{(\gamma)\prime}=\gamma W^{(\gamma)}$, proving \eqref{eq:f3-piecewise}.  The cases $d\le0$ and $d\ge r$ follow from the range restriction of the post-supremum path.
    
(iii) Proof was given in \cite{VardarAcar2021}.
%\end{proof}

%\section{Brownian consistency check}\label{sec:brownian}
\begin{remark}[Comparison of the pre-$H_I$ process 
with (\cite{Salminen2007}, Thm 3.5 and Prop 4.1)]
Let $X$ be standard Brownian motion and let $T\sim\operatorname{Exp}(\gamma)$.  We have
\begin{equation*}\label{eq:brownian-scale}
W^{(\gamma)}(x)=\frac{2}{\sqrt{2\gamma}}\sinh(\sqrt{2\gamma} x),
\qquad
Z^{(\gamma)}(x)=\cosh(\sqrt{2\gamma} x),
\qquad x\ge0.
\end{equation*}
Substitution into the $h$-functions above recovers the Brownian decompositions of \cite{Salminen2007}. 

\begin{eqnarray}\label{eq:cek1}
h_1(x)&=&
\frac{W^{(\gamma)\prime}(x-a)W^{(\gamma)}(b-a)-W^{(\gamma)\prime}(b-a)W^{(\gamma)}(x-a)}
     {\big(W^{(\gamma)}(b-a)\big)^2}\\\nonumber&=&\frac{\sqrt{2\gamma}\sinh(\sqrt{2\gamma} (x-b))}{\sinh^2(\sqrt{2\gamma} (b-a))},\quad a<x<b\nonumber
\end{eqnarray}
It follows that
\begin{eqnarray}\label{eq:cek2}
\frac{h'_1(x)}{h_1(x)}&=&\frac{\sqrt{2\gamma}\cosh(\sqrt{2\gamma} (x-b))}{\sinh(\sqrt{2\gamma} (x-b))}=\frac{\sqrt{2\gamma}(1+e^{-2\sqrt{2\gamma}(x-b)})}{1+e^{-2\sqrt{2\gamma}(x-b)}},
\quad a<x<b\nonumber
\end{eqnarray}
and we have $\frac{\gamma h_1(x)-\frac{1}{2}h''_1(x)}{h_1(x)}=0;$
therefore, the generator of the $h$-transform coincides with the generator Brownian motion with drift $\sqrt{2\gamma}$ started from 0 conditioned to hit $a$ before $b$ and killed when it hits $a$ as found in \cite{Salminen2007}.
%For example, \eqref{eq:h-pre-infimum} becomes, up to an irrelevant multiplicative constant,
%\[
%h_0(x)=e^{-\sqrt{2\gamma}(x-a)},
%\]
%which transforms the Brownian generator $\frac12 f''$ into $\frac12 f''-\sqrt{2\gamma} f'$, i.e. Brownian motion with drift $-\sqrt{2\gamma}$ killed at the lower boundary.
The final distribution formulae also reduce to explicit Brownian expressions.  With $r=b-a$, \ref{eq:f1-piecewise} gives, for $-a<d<r$,
\begin{equation}\label{eq:brownian-f1}
\Pp_{a,b}^{<}\{M_{0,H_I}^-<d\}
=
\frac{\sinh(\sqrt{2\gamma} r)}{\sinh(\sqrt{2\gamma} b)}
\frac{\sinh(\sqrt{2\gamma}(d+a))}{\sinh(\sqrt{2\gamma} d)}.
\end{equation}
Similarly, \ref{eq:f2-piecewise} gives, for $0<d<r$,
\begin{equation}\label{eq:brownian-f2}
\Pp_{a,b}^{<}\{M_{H_I,H_S}^-<d\}
=
\frac{\sinh(\sqrt{2\gamma} r)}{\sinh(\sqrt{2\gamma} d)}
\exp\{-\sqrt{2\gamma}(r-d)\coth(\sqrt{2\gamma} d)\},
\end{equation}
and \ref{eq:f3-piecewise} gives
\begin{equation}\label{eq:brownian-f3}
\Pp_{a,b}^{<}\{M_{H_S,T}^-<d\}
=
\frac{\tanh(\sqrt{2\gamma} d/2)}{\tanh(\sqrt{2\gamma} r/2)}.
\end{equation}
These expressions have the correct boundary limits at $d=0$, $d=-a$ and $d=r$, and match the corresponding Brownian path-decomposition formulae \cite{Salminen2007}. 
\end{remark}
\section{Conclusion}\label{sec:conclusion}

We have given a self-contained scale-function formulation of the conditional path decomposition needed to compute maximum drawdowns under the ordering $H_I<H_S$ for spectrally negative L\'evy processes.  The main theorem proves the conditional maximum-drawdown law of the pre-infimum process and identifies the laws of intermediate and post-supremum  processes.   We have verified that the resulting expressions define legitimate distribution functions with consistent boundary limits. Our contributions for spectrally negative L\'evy processes confirm the explicit formulae found in \cite{Salminen2007} for the special case of Brownian motion.

\section*{Acknowledgements}
% The authors thank the referees for helpful comments that improved the presentation.  Funding and conflict-of-interest statements may be added before submission, if applicable.
This work is supported by Tübitak Project with number 124F094.

%\begin{thebibliography}{99}
\bibliography{myBiblio}

@article{Avram2019,
  title={The $w$,$z/\nu$,$\delta$ paradigm for the first passage of strong Markov processes without positive jumps},
  author={Avram, Florin and Grahovac, Danijel and Vardar-Acar, Ceren},
  journal={Risks},
  volume={7},
  number={1},
  year={2019},
  doi={10.3390/risks7010010}
}

@article{Avram2020,
  title={The $W$, $Z$ scale functions kit for first passage problems of spectrally negative Lévy processes, and applications to control problems},
  author={Avram, Florin and Grahovac, Danijel and Vardar-Acar, Ceren},
  journal={ESAIM: Probability and Statistics},
  volume={24},
  pages={454--525},
  year={2020},
  doi={10.1051/ps/2019025}
}

@article{Avram2004,
  title={Exit problems for spectrally negative Lévy processes and applications to (Canadized) Russian options},
  author={Avram, Florin and Kyprianou, Andreas E. and Pistorius, Martijn R.},
  journal={The Annals of Applied Probability},
  volume={14},
  number={1},
  pages={215--238},
  year={2004},
  doi={10.1214/aoap/1075828052}
}

@article{Baurdoux2017,
  title={On future drawdowns of Lévy processes},
  author={Baurdoux, Erik J. and Palmowski, Zbigniew and Pistorius, Martijn R.},
  journal={Stochastic Processes and their Applications},
  volume={127},
  number={8},
  pages={2679--2698},
  year={2017},
  doi={10.1016/j.spa.2016.12.008}
}

@book{Bertoin1996,
  title={Lévy Processes},
  author={Bertoin, Jean},
  year={1996},
  publisher={Cambridge University Press},
  address={Cambridge},
  series={Cambridge Tracts in Mathematics},
  volume={121},
  isbn={978-0-521-64632-1}
}

@article{Carr2011,
  title={Maximum drawdown insurance},
  author={Carr, Peter and Zhang, Hailiang and Hadjiliadis, Olympia},
  journal={International Journal of Theoretical and Applied Finance},
  volume={14},
  number={8},
  pages={1195--1230},
  year={2011},
  doi={10.1142/S0219024911006683}
}

@article{Chekhlov2005,
  title={Drawdown measure in portfolio optimization},
  author={Chekhlov, Alexei and Uryasev, Stanislav and Zabarankin, Michael},
  journal={International Journal of Theoretical and Applied Finance},
  volume={8},
  number={1},
  pages={13--58},
  year={2005},
  doi={10.1142/S0219024905002767}
}

@article{Douady2000,
  title={On probability characteristics of downfalls in a standard Brownian motion},
  author={Douady, R. and Shiryaev, A. N. and Yor, M.},
  journal={Theory of Probability \& Its Applications},
  volume={44},
  number={1},
  pages={29--38},
  year={2000},
  doi={10.1137/S0040585X97978208}
}

@article{Hadjiliadis2006,
  title={Drawdowns preceding rallies in the Brownian motion model},
  author={Hadjiliadis, Olympia and Večeř, Jan},
  journal={Quantitative Finance},
  volume={6},
  number={5},
  pages={403--409},
  year={2006},
  doi={10.1080/14697680600764227}
}

@article{Huang2013,
  title={Multivariate risk models under heavy-tailed risks},
  author={Huang, Wei and Weng, Cheng and Zhang, Yichuan},
  journal={Applied Stochastic Models in Business and Industry},
  year={2013},
  doi={10.1002/asmb.1987}
}

@incollection{Kuznetsov2013,
  title={The theory of scale functions for spectrally negative Lévy processes},
  author={Kuznetsov, Alexey and Kyprianou, Andreas E. and Rivero, Victor},
  booktitle={Lévy Matters II},
  pages={97--186},
  year={2013},
  publisher={Springer},
  doi={10.1007/978-3-642-31407-0_2}
}

@book{Kyprianou2014,
  title={Fluctuations of Lévy Processes with Applications: Introductory Lectures},
  author={Kyprianou, Andreas E.},
  year={2014},
  publisher={Springer Science \& Business Media},
  doi={10.1007/978-3-642-37632-0}
}

@article{Landriault2015,
  title={Analysis of a drawdown-based regime-switching Lévy insurance model},
  author={Landriault, David and Li, Bin and Li, Shuanming},
  journal={Insurance: Mathematics and Economics},
  volume={60},
  pages={98--107},
  year={2015},
  doi={10.1016/j.insmatheco.2014.10.005}
}

@article{Leal2005,
  title={Maximum drawdown},
  author={Leal, Ricardo P. C. and Mendes, Beatriz V. de M.},
  journal={The Journal of Alternative Investments},
  volume={7},
  number={4},
  pages={83--91},
  year={2005},
  doi={10.3905/jai.2005.491503}
}

@article{Magdon-Ismail2004,
  title={On the maximum drawdown of a Brownian motion},
  author={Magdon-Ismail, Malik and Atiya, Amir F. and Pratap, Anjali and Abu-Mostafa, Yaser S.},
  journal={Journal of Applied Probability},
  volume={41},
  number={1},
  pages={147--161},
  year={2004},
  doi={10.1239/jap/1077134835}
}

@article{Mijatovic2012,
  title={On the drawdown of completely asymmetric Lévy processes},
  author={Mijatović, Aleksandar and Pistorius, Martijn R.},
  journal={Stochastic Processes and their Applications},
  volume={122},
  number={11},
  pages={3812--3836},
  year={2012},
  doi={10.1016/j.spa.2012.07.011}
}

@article{Millar1977,
  title={Zero-one laws and the minimum of a Markov process},
  author={Millar, P. W.},
  journal={Transactions of the American Mathematical Society},
  volume={226},
  pages={365--391},
  year={1977},
  publisher={American Mathematical Society},
  doi={10.1090/S0002-9947-1977-0433606-6}
}

@article{Pospisil2010,
  title={Portfolio sensitivity to changes in the maximum and the maximum drawdown},
  author={Pospisil, Libor and Vecer, Jan},
  journal={Quantitative Finance},
  volume={10},
  number={6},
  pages={617--627},
  year={2010},
  publisher={Taylor \& Francis},
  doi={10.1080/14697680903008751}
}

@article{Rossello2021,
  title={A refined measure of conditional maximum drawdown},
  author={Rossello, David and Lo Cascio, Salvatore},
  journal={Risk Management},
  volume={23},
  number={4},
  pages={301--321},
  year={2021},
  publisher={Palgrave Macmillan},
  doi={10.1057/s41283-021-00074-0}
}

@article{Salminen2007,
  title={On maximum increase and decrease of Brownian motion},
  author={Salminen, Paavo and Vallois, Pierre},
  journal={Annales de l'Institut Henri Poincaré (B) Probability and Statistics},
  volume={43},
  pages={655--676},
  year={2007},
  publisher={Elsevier},
  doi={10.1016/j.anihpb.2006.09.001}
}

@book{Sornette2003,
  title={Why Stock Markets Crash: Critical Events in Complex Financial Systems},
  author={Sornette, Didier},
  year={2003},
  publisher={Princeton University Press},
  address={Princeton, NJ}
}

@article{VardarAcar2017,
  title={Maximum loss and maximum gain of spectrally negative Lévy processes},
  author={Vardar-Acar, Ceren and Çağlar, Mahir},
  journal={Extremes},
  volume={20},
  number={2},
  pages={301--308},
  year={2017},
  publisher={Springer},
  doi={10.1007/s10687-016-0264-1}
}

@article{VardarAcar2021,
  title={Maximum drawdown and drawdown duration of spectrally negative Lévy processes decomposed at extremes},
  author={Vardar-Acar, Ceren and Çağlar, Mahir and Avram, Florin},
  journal={Journal of Theoretical Probability},
  volume={34},
  pages={1486--1505},
  year={2021},
  publisher={Springer},
  doi={10.1007/s10959-020-01034-0}
}

@inproceedings{Vecer2006,
  title={Maximum drawdown and directional trading},
  author={Vecer, Jan},
  booktitle={Frankfurt Math Finance Workshop Paper},
  year={2006}
}

@article{Vecer2007,
  title={Preventing portfolio losses by hedging maximum drawdown},
  author={Vecer, Jan},
  journal={Wilmott},
  volume={5},
  number={4},
  pages={1--8},
  year={2007}
}

@article{Zhang2010,
  title={Drawdowns and rallies in a finite time-horizon},
  author={Zhang, Hongzhong and Hadjiliadis, Olympia},
  journal={Methodology and Computing in Applied Probability},
  volume={12},
  number={2},
  pages={293--308},
  year={2010},
  publisher={Springer},
  doi={10.1007/s11009-009-9135-4}
}

@article{grossman1993optimal,
  title={Optimal investment strategies for controlling drawdowns},
  author={Grossman, Sanford J and Zhou, Zhongquan},
  journal={Mathematical finance},
  volume={3},
  number={3},
  pages={241--276},
  year={1993},
  publisher={Wiley Online Library}
}

@article{cvitanic1995portfolio,
  title={On portfolio optimization under "drawdown" constraints},
  author={Cvitanic, Jaksa and Karatzas, Ioannis},
  journal={IMA Volume in Mathematical Finance},
  volume={65},
  pages={35--42},
  year={1995}
}

@article{SALMINEN20205592,
title = {On the maximum increase and decrease of one-dimensional diffusions},
journal = {Stochastic Processes and their Applications},
volume = {130},
number = {9},
pages = {5592-5604},
year = {2020},
issn = {0304-4149},
doi = {https://doi.org/10.1016/j.spa.2020.04.001}
}

@article{salminen2025drawdowns,
  title={Drawdowns of diffusions},
  author={Salminen, Paavo and Vallois, Pierre},
  journal={ESAIM: Probability and Statistics},
  volume={29},
  pages={357--380},
  year={2025},
  publisher={EDP Sciences}
}

@article{ZHANG2023104669,
title = {A general method for analysis and valuation of drawdown risk},
journal = {Journal of Economic Dynamics and Control},
volume = {152},
pages = {104669},
year = {2023},
issn = {0165-1889},
doi = {https://doi.org/10.1016/j.jedc.2023.104669}
}

@article{10.1214/ECP.v20-3945,
author = {Yueyun Hu and Zhan Shi and Marc Yor},
title = {{The maximal drawdown of the Brownian meander}},
volume = {20},
journal = {Electronic Communications in Probability},
number = {none},
publisher = {Institute of Mathematical Statistics and Bernoulli Society},
pages = {1 -- 6},
year = {2015},
doi = {10.1214/ECP.v20-3945}
}

@Article{risks7040105,
AUTHOR = {Mayerhofer, Eberhard},
TITLE = {Three Essays on Stopping},
JOURNAL = {Risks},
VOLUME = {7},
YEAR = {2019},
NUMBER = {4},
ARTICLE-NUMBER = {105}
}
\bibliographystyle{abbrvnat}

\end{document}